\documentclass{gSTA2e}

\usepackage{latexsym}\usepackage{amssymb,latexsym}
\usepackage[mathscr]{eucal} \usepackage{amsthm}
\usepackage{amsmath}
\numberwithin{equation}{section}
\numberwithin{table}{section}
\usepackage{amsfonts} \usepackage{layout}
\usepackage{bold-extra}
\usepackage{multirow}
\newtheorem{theorem}{Theorem}

\usepackage{mathrsfs}
\usepackage{xcolor}
\usepackage{listings}
\usepackage{subfigure}
\usepackage{graphicx}
\begin{document}

\title{\textit{Construction of Simultaneous Confidence Bands for Multiple Logistic Regression
  Models over Restricted Regions}}

\author{
\name{Lucy Kerns$^\ast$ \thanks{$^\ast$ Email: xlu@ysu.edu}}
\affil{Department of Mathematics and Statistics, \\Youngstown State University, Youngstown OH 44555}
}
\maketitle

\begin{abstract}
This article presents methods for constructing an asymptotic
hyperbolic band under the multiple logistic
regression model when the predictor variables are restricted
to a specific region~$\mathscr{X}$.  Scheff\'{e}'s method yields unnecessarily wide,  and hence conservative,  bands if the predictor variables can be restricted to a certain region.
Piegorsch and Casella \cite{PC88} developed a procedure to build an asymptotic confidence band for the multiple
logistic regression model over particular regions. Those regions are shown to be special cases of the region~$\mathscr{X}$, which was first investigated by Seppanen and Uusipaikka \cite{SU92} in the multiple linear regression context.  This article also provides methods for constructing conservative confidence bands when the restricted region is not of the specified form. Particularly, rectangular restricted regions, which are commonly encountered in practice, are considered.  Two examples are given to illustrate the proposed methodology, and one example shows that the proposed procedure outperforms the method given by Piegorsch and Casella \cite{PC88}. \bigskip
\end{abstract}


\begin{classcode}
 62J12; 62J05
\end{classcode}

\section {Introduction}
The logistic regression model is a statistical method for analyzing the effects of independent variables on a dichotomous dependent response, that is, a
response that takes values 1 (success, normal, positive, etc.) and 0 (failure, abnormal, negative, etc.). The logistic model specifies the response probability as
\begin{equation}\label{eq:1.1}
P(Y = 1) = p(\mathbf{x})=  1/[1 + \text{exp}(-\mathbf{x}'\boldsymbol{\beta})],
\end{equation}
where $Y$ is the response variable, $\mathbf{x} = (1,\,\, x_1, \ldots, x_{p-1})'$ with $x_1, x_2,\ldots, x_{p-1}$ being the set of predictor variables, and
$\boldsymbol{\beta}=(\beta_0,\,\, \beta_1, \dots, \beta_{p-1})'$ is the $p \times 1$ vector of unknown parameters.

The parameter $\boldsymbol{\beta}$ in the logistic models can be estimated using the method of maximum likelihood.  The ML estimation of $\boldsymbol{\beta}$ requires iterative computational methods, but existing computer software has facilitated the computation. The ML estimator, denoted by $\hat{\boldsymbol{\beta}}$,  under certain
regularity conditions (Kendall and Stuart \cite{KS79}),
 follows asymptotically as:
\[
\hat{\boldsymbol{\beta}} \sim ASN_p(\boldsymbol{\beta},\boldsymbol{F^{-1}}),
\]
where $\boldsymbol{F}^{-1}$ is the inverse of the Fisher information matrix, which is commonly provided by computer software.

It is often of interest to construct a simultaneous confidence band on $p(\mathbf{x})$  because it provides useful information on the plausible range of the unknown response probability. The  construction of the band can be simplified by using the logit link function
\[
\text{logit}(p(\mathbf{x}))= \text{log}_e\left[\frac{p(\mathbf{x})}{1-p(\mathbf{x})}\right]= \mathbf{x}'\boldsymbol{\beta}.
\]
Applying the logit function on $p(\mathbf{x})$ enables us to  transform the problem
of constructing confidence bands for the response probability $p(\mathbf{x})$
to the problem of constructing confidence bands for the linear
predictor $\mathbf{x}'\boldsymbol{\beta}$ on which the bands are
defined.  A $100(1-\alpha)\%$ two-sided hyperbolic band for the linear predictor
$\mathbf{x}'\boldsymbol{\beta}$ has the form
\begin{equation}\label{eq:1.2}
\mathbf{x}'\boldsymbol{\beta} \in \mathbf{x}'\boldsymbol{\hat{\beta}} \pm c(\mathbf{x}'\boldsymbol{F}^{-1}\mathbf{x})^{1/2},
\end{equation}
where $c$ satisfies
\begin{equation}\label{eq:1.3}
P[\mathbf{x}'\boldsymbol{\beta} \in \mathbf{x}'\boldsymbol{\hat{\beta}} \pm c(\mathbf{x}'\boldsymbol{F}^{-1}\mathbf{x})^{1/2}] = 1 - \alpha.
\end{equation}

Since the logit function is a monotonically increasing
function, a $100(1-\alpha)\%$ two-sided confidence band for $p(\mathbf{x})$ in
the logistic regression model is then given by
\begin{equation}\label{eq:1.4}
 \{1 + \text{exp}[-\mathbf{x}'\boldsymbol{\hat{\beta}} + c(\mathbf{x}'\boldsymbol{F}^{-1}\mathbf{x})^{1/2}]\} ^{-1} \leq p(\mathbf{x}) \leq
 \{1 + \text{exp}[-\mathbf{x}'\boldsymbol{\hat{\beta}} - c(\mathbf{x}'\boldsymbol{F}^{-1}\mathbf{x})^{1/2}]\} ^{-1}.
\end{equation}

The work of constructing confidence bands in simple and multiple linear regression models dates back to Working $\&$ Hotelling \cite{WH29}. The confidence band developed by Scheff\'{e}~\cite{S59} is the well known two-sided hyperbolic band when there are no restrictions on the values of the predictor variables.  Scheff\'{e}'s bands can be unnecessarily wide,  and hence conservative,  if it is reasonable for us to constrain the predictor variables to a certain region. A considerable amount of work has been done on improving the earlier work by restricting the predictor variables to some specified regions. For example, Halperin $\&$ Gurain \cite{HG68} and Liu $\&$ Lin \cite{LL09} constructed two-sided hyperbolic confidence bands over an ellipsoidal region.
Wynn $\&$ Bloomfield \cite{WB71} and  Uusipaikka \cite{U83} developed exact two-sided hyperbolic bands when the predictor
variable is restricted to an interval or union of intervals. Casella and
Strawderman \cite{CS80} were able to build an exact two-sided band over a region that was more general than the region considered by Halperin $\&$ Gurain \cite{HG68} and Liu $\&$ Lin \cite{LL09}.
Seppanen and Uusipaikka \cite{SU92} further studied the region investigated by Casella and
Strawderman \cite{CS80}, and  developed an exact two-sided confidence band over a region that is even more general than the one considered by Casella and
Strawderman. Liu \emph{et al.} \cite{LJZD05}
provided simulation-based confidence bands for  multiple regression models when certain restrictions are placed on the predictor variables.  Liu's recent book \cite{L10} gave a comprehensive overview of the methodology for constructing simultaneous confidence bands and the applications of the bands in various statistical problems.

Much less work has been done on constructing confidence bands for the
logistic regression model. Brand, Pinnock, and Jackson \cite{BPJ73}
constructed confidence bands for both $p(\mathbf{x})$ and the inverse of $p(\mathbf{x})$
in the simple logistic case with no restrictions on the predictor variable. Hauck \cite{H83} extended their work to more than one
predictor variable yet still with no constraints. When constraints exist on the predictor variables, Piegorsch and Casella \cite{PC88} developed asymptotic two-sided bands by extending
the early work of Casella and Strawderman \cite{CS80} from the multiple linear
regression model to the multiple logistic regression model. Wei Liu's book \cite[Chapter 8]{L10} presented confidence bands for the logistic model
with more than one explanatory variable. The method in the book utilized
simulation-based confidence bands (Liu \emph{et al.} \cite{LJZD05}) for the
linear predictor $\mathbf{x}'\boldsymbol{\beta}$ in the multiple
linear regression model, and the desired bands for the logistic model were then
obtained via the logit link function.  Kerns \cite{K15} considered the simple logistic regression case, and was able to  develop  asymptotic
two-sided and one-sided simultaneous hyperbolic bands when the predictor variable is restricted
to a given interval, such as ($\emph{l},\,\, \emph{u}$) with $l$ and $u$ being given real numbers.

The constraint region considered  by Piegorsch and Casella \cite{PC88} is the same region studied by Casella and Strawderman \cite{CS80}. In this paper, we will focus on building a $100(1-\alpha)\%$ confidence band over the  region that was previous studied by Seppanen and Uusipaikka~\cite{SU92}, and as mentioned earlier, this region is more general than the one discussed by Casella $\&$ Strawderman \cite{CS80} and Piegorsch $\&$ Casella \cite{PC88}. The constraint region has the form
\begin{equation}\label{eq:1.5}
\mathscr{X} = \{\mathbf{x}: \mathbf{x}'\boldsymbol{F}^{-1}\boldsymbol{Z}(\boldsymbol{Z}'\boldsymbol{F}^{-1}\boldsymbol{Z})^{-1}
\boldsymbol{Z}'\boldsymbol{F}^{-1}\mathbf{x} \geq a^2\mathbf{x}'\boldsymbol{F}^{-1}\mathbf{x}\},
\end{equation}
where $\boldsymbol{Z} = (\boldsymbol{z}_1, \boldsymbol{z}_2, \ldots, \boldsymbol{z}_r)$ is an arbitrary given $p \times r$ matrix, whose columns are linearly independent, and $a$ is a given real number such that $0 \leq a \leq 1$. This region can be further written as
\begin{equation}\label{eq:1.6}
\mathscr{X} = \{\mathbf{x}: \rho(\mathbf{x},E) \geq a\},
\end{equation}
where
\begin{equation}\label{eq:1.7}
\rho(\mathbf{x},E)= \Big\{\frac{ \mathbf{x}'\boldsymbol{F}^{-1}\boldsymbol{Z}(\boldsymbol{Z}'\boldsymbol{F}^{-1}\boldsymbol{Z})^{-1}
\boldsymbol{Z}'\boldsymbol{F}^{-1}\mathbf{x}}{\mathbf{x}'\boldsymbol{F}^{-1}\mathbf{x}}\Big\}^{1/2}
\end{equation}
is the multiple correlation coefficient between the random variable $\mathbf{x}'\boldsymbol{\hat{\beta}}$ and the random vector $\boldsymbol{Z}'\boldsymbol{\hat{\beta}}$, and $E$ is a $r$-dimensional subspace of $\mathbb{R}^p$ ($1 \leq r \leq p$), which is spanned by the columns of $\boldsymbol{Z}$.

In what follows, we will focus on the region $\mathscr{X}$ defined in Equation~(\ref{eq:1.5}), or equivalently Equation~(\ref{eq:1.6}),  and propose a method for calculating the critical value $c_{\alpha}$ in a $100(1-\alpha)\%$ two-sided band:
\begin{equation}\label{eq:1.8}
\mathbf{x}'\boldsymbol{\beta} \in \mathbf{x}'\boldsymbol{\hat{\beta}} \pm c_{\alpha}(\mathbf{x}'\boldsymbol{F}^{-1}\mathbf{x})^{1/2},
\end{equation}
where $c_{\alpha}$ satisfies
\begin{equation}\label{eq:1.9}
P[\mathbf{x}'\boldsymbol{\beta} \in \mathbf{x}'\boldsymbol{\hat{\beta}} \pm c_{\alpha}(\mathbf{x}'\boldsymbol{F}^{-1}\mathbf{x})^{1/2},\quad \forall \,\, \mathbf{x} \in \mathscr{X} ] = 1 - \alpha.
\end{equation}
Then, a $100(1-\alpha)\%$ two-sided confidence band for $p(\mathbf{x})$ in
the multiple logistic regression model when the predictor variables are constrained to the given region $\mathscr{X}$ can be obtained as
\begin{equation}\label{eq:1.10}
 \{1 + \text{exp}[-\mathbf{x}'\boldsymbol{\hat{\beta}} + c_{\alpha}(\mathbf{x}'\boldsymbol{F}^{-1}\mathbf{x})^{1/2}]\} ^{-1} \leq p(\mathbf{x}) \leq
 \{1 + \text{exp}[-\mathbf{x}'\boldsymbol{\hat{\beta}} - c_{\alpha}(\mathbf{x}'\boldsymbol{F}^{-1}\mathbf{x})^{1/2}]\} ^{-1}.
\end{equation}

In this paper, we provide explicit expressions for determining the critical values in the multiple logistic regression setting when the predictor variables are constrained to the region $\mathscr{X}$. Compared with the region studied by Piegorsch and Casella \cite{PC88}, our region is more general, and hence covers more different forms of constraints. Consequently, our proposed method is able to produce narrower confidence bands if the restricted region is of the form $\mathscr{X}$, but not of the form considered by Piegorsch and Casella. In particular, we will show that the interval restriction on the predictor variable in the simple logistic regression model belongs to the class of the region $\mathscr{X}$. As Piegorsch and Casella noted, however, the regions of this form cannot be recovered from their regions. Instead, they had to apply the embedding procedure from Casella and Strawderman \cite{CS80}. Not only did the embedding procedure require a fair amount of computation, but it also produced wider bands. Hence, in the simple logistic regression setting, our proposed method has an advantage over theirs.  Wei Liu's method presented in his book \cite[Chapter 8]{L10} is very broad, but relies on simulation. Our method is more focused but admits tractable forms.

The logistic regression model belongs to the family of the
generalized linear model (GLM), and the methodology proposed here can
also be applied to other forms of GLM (the probit model and the complementary-log-log
model, for example), that can be transformed into the standard
regression model via a link function. In each case, a set of
simulations are required to confirm the validity of the method in
small samples, this would drive the size of this paper to unwieldy levels. Thus in the interest of brevity, we devote our attention to
the logistic model to illustrate the methodology.

This paper is organized as follows. We present the results on calculating the critical value in Section 2 when the predictor variables are restricted to the region~(\ref{eq:1.5}). We will also discuss how to obtain conservative confidence bands when the constrained region is not of the form~(\ref{eq:1.5}). Two examples are given in Section 3 to illustrate the proposed methodology. In Section 4, a Monte Carlo simulation is run to
investigate how well the asymptotic approximation holds for small
sample sizes.

\section{Theory and Methods} \label{theory}
\subsection{Two-sided Bands Over Constraint Regions}
Theorem~1 in Seppanen $\&$ Uusipaikka \cite{SU92} provided a method for constructing exact bands when the region of interest is of the form~(\ref{eq:1.5}) in the multiple linear regression model. We will modify their methodology and extend their work to fit the multiple logistic case. Theorem~1 in this paper presents an explicit expression for calculating the critical value $c_{\alpha}$ in the multiple logistic model with the same constraint region.  Since our methodology is based on the one proposed by Seppanen $\&$ Uusipaikka,  there are similarities between the theorems and the proofs given here and in their paper. The main difference originates from the definition of the vectors $\boldsymbol{u}$ and $\boldsymbol{v}$ (see the proof below). These two vectors were jointly distributed as $\boldsymbol{N}_p(\boldsymbol{0},\sigma^2\boldsymbol{I}_p)$ in the multiple linear regression model, while in the multiple logistic setting, their joint asymptotic distribution is $\boldsymbol{N}_p(\boldsymbol{0},\boldsymbol{I}_p)$. As a result, the random variable $w$ defined by $w = \|\boldsymbol{u}\|^2 + \|\boldsymbol{v}\|^2$  follows a $\chi^2$ distribution with $p$ degrees of freedom, while $w$, which was defined slightly differently in their paper, followed an $F$ distribution with $p$ and $n-p$ degrees of freedom.  Consequently,  the  chi-square distribution, rather than the $F$ distribution, is involved in our final results.  In the interest of completeness, we include most of the computational details regarding the critical value $c_{\alpha}$.

\begin{theorem}
The two-sided confidence bands~(\ref{eq:1.8}) and (\ref{eq:1.10}),  when the predictor variables are constrained to the region $\mathscr{X}$,  are large-sample $100(1 - \alpha)\%$ confidence bands if the critical value $c_{\alpha}$ satisfies $P(G \leq c_{\alpha}^2) = 1 - \alpha$, where the random variable $G$ has the distribution function
\begin{equation}\label{eq:2.1}
P(G \leq g) = F(g) + \int_g^{g/(1-a^2)}H(m(\sqrt{g/w}))f(w)dw,
\end{equation}
where $F$ and $f$ are the distribution and density functions of the $\chi^2$ distribution with $p$ degrees of freedom, $H$ is the distribution function of the beta($r/2, (p-r)/2$) distribution, and $m$ is the function
\[
m(t) =  \{at - [(1-a^2)(1-t^2)]^{1/2}\}^2, \,\,\, \forall  \,0 \leq t \leq 1.
\]
\end{theorem}

\emph{\textbf{Remark 1}}. It is clear that the region $\mathscr{X}$ is invariant if $\boldsymbol{Z}$ is replaced by $\boldsymbol{ZA}$ for any non-singular $r \times r$ matrix $\boldsymbol{A}$, therefore it depends only on the subspace spanned by the columns of $\boldsymbol{Z}$.  Seppanen $\&$ Uusipaikka \cite{SU92} stated, without proof, that the constraints considered by Casella $\&$ Strawderman \cite{CS80} is a special case of the region $\mathscr{X}$, in which the  columns of $\boldsymbol{Z}$ consist of $r$ orthonormal eigenvectors of $\boldsymbol{F}^{-1}$. We will explain this claim in more details here.

The inverse of the Fisher information matrix $\boldsymbol{F}^{-1}$ is positive semi-definite, and can be written as $\boldsymbol{F}^{-1} = \boldsymbol{PDP}'$, where $\boldsymbol{D}$ = diag$\{\lambda_i\}$ is the diagonal matrix of the eigenvalues of $\boldsymbol{F}^{-1}$, and $\boldsymbol{P}$ is the matrix of corresponding orthonormal eigenvectors.

 If we choose $\boldsymbol{Z}$ in~(\ref{eq:1.5}) to be $\boldsymbol{PD}^{-1/2}$, then $(\boldsymbol{P}\boldsymbol{D}^{-1/2})'\boldsymbol{F}^{-1}(\boldsymbol{P}\boldsymbol{D}^{-1/2}) = \boldsymbol{I}_r$, and  there exists a $p \times (p-r)$ matrix $\boldsymbol{M}$ such that $\boldsymbol{M}'\boldsymbol{F}^{-1}\boldsymbol{M} = \boldsymbol{I}_{p-r}$ and $(\boldsymbol{P}\boldsymbol{D}^{-1/2})'\boldsymbol{F}^{-1}\boldsymbol{M} = \boldsymbol{0}_{r \times (p-r)}$.

Denote $\boldsymbol{z} = (\boldsymbol{P}\boldsymbol{D}^{-1/2})'\boldsymbol{F}^{-1}\mathbf{x}$ and $\boldsymbol{d} = \boldsymbol{M}'\boldsymbol{F}^{-1}\mathbf{x}$. Then
\[
\mathbf{x}'\boldsymbol{F}^{-1}(\boldsymbol{P}\boldsymbol{D}^{-1/2})[(\boldsymbol{P}\boldsymbol{D}^{-1/2})'\boldsymbol{F}^{-1}(\boldsymbol{Z}\boldsymbol{P}^{-1/2})]^{-1}
(\boldsymbol{P}\boldsymbol{D}^{-1/2})'\boldsymbol{F}^{-1}\mathbf{x} = \|\boldsymbol{z}\|^2,
\]
and
\[
\mathbf{x}'\boldsymbol{F}^{-1}\mathbf{x} =  \|\boldsymbol{z}\|^2 +  \|\boldsymbol{d}\|^2.
\]
That is, the region $\mathscr{X}$ can be written as $\{(\boldsymbol{z},\boldsymbol{d} )': \|\boldsymbol{z}\|^2 / (\|\boldsymbol{z}\|^2 +  \|\boldsymbol{d}\|^2) \geq a^2\}$, or equivalently, $\{(\boldsymbol{z},\boldsymbol{d} )': \|\boldsymbol{z}\|^2  \geq \frac{a^2}{1 - a^2} \|\boldsymbol{d}\|^2\}$. This is the same region that Casella $\&$ Strawderman~\cite{CS80} and Piegorsch $\&$ Casella \cite{PC88} considered in their articles.

\vspace{.3cm}
\emph{\textbf{Remark 2}}. If the subspace $E$ is one-dimensional ($r = 1$), that is, if $E$ = span\{$\boldsymbol{z}$\}, where $\boldsymbol{z}$ is a $p \times 1$ vector, then the region $\mathscr{X}$ can be written as:
\[
\mathscr{X} = \{\mathbf{x}: a^2\mathbf{x}'\boldsymbol{F}^{-1}\mathbf{x}\boldsymbol{z}'\boldsymbol{F}^{-1}\boldsymbol{z} \leq (\mathbf{x}'\boldsymbol{F}^{-1}\boldsymbol{z})^2\} =
\{\mathbf{x}:|\rho(\mathbf{x}, \boldsymbol{z})| \geq a\},
\]
where $\rho(\mathbf{x}, \boldsymbol{z})$ is the correlation coefficient between the random variables $\mathbf{x}'\hat{\boldsymbol{\beta}}$ and $\boldsymbol{z}'\hat{\boldsymbol{\beta}}$.

For example, consider the simple logistic regression model:
\[p(\mathbf{x})= 1/[1 + \text{exp}(-\mathbf{x}'\boldsymbol{\beta})] = 1/[1 + \text{exp}(-(\beta_0 + \beta_1x))],\]
which can be written in the more general form
\[p(\mathbf{x})= 1/[1 + \text{exp}(-(\beta_0x_0 + \beta_1x_1))].\]
If the predictor variable $x$ is restricted to lie in an interval, that is, if $x \in (\emph{l},\,\emph{u})$, where $\emph{l}$ and $\emph{u}$ are pre-specified constants, then the restriction over the interval is equivalent to restricting $(x_0,\,\, x_1)'$ to be in the set, $\{(x_0,\,\,x_1): x_0 = 1, \,\, \emph{l} < x_1 < \emph{u}\}$. This set is a line segment, as illustrated in Figure~1.

\begin{figure}
\begin{center}
\subfigure[Figure 1]{\resizebox*{7cm}{!}{\includegraphics{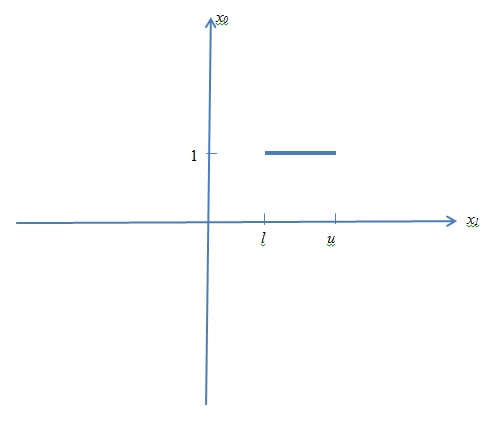}}}\hspace{15pt}
\subfigure[Figure 2]{\resizebox*{7cm}{!}{\includegraphics{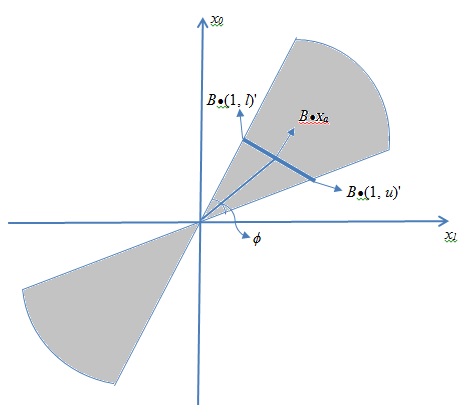}}}
\end{center}
\end{figure}

The matrix $\boldsymbol{F^{-1}}$ can also be written as $\boldsymbol{F^{-1}} = \boldsymbol{B}^2$, where $\boldsymbol{B} = \boldsymbol{P}\boldsymbol{D}^{1/2}\boldsymbol{P}'$. Then, as Kerns \cite{K15} pointed out, for the two-sided band in Equation~(\ref{eq:1.8}), the simultaneous confidence level is given by

\begin{align}\label{eq:add}
 &P[\mathbf{x}'\boldsymbol{\beta} \in \mathbf{x}'\boldsymbol{\hat{\beta}} \pm c_{\alpha}(\mathbf{x}'\boldsymbol{F}^{-1}\mathbf{x})^{1/2}, \,\,\text{for all}\,\, x \in (\emph{l},\emph{u})] \notag\\
 &=P\left[\sup\limits_{x \in (\emph{l},\emph{u})}
     \frac{|\mathbf{x}'\boldsymbol{\beta} - \mathbf{x}'\boldsymbol{\hat{\beta}}|}{(\mathbf{x}'\boldsymbol{F}^{-1}\mathbf{x})^{1/2}} < c_{\alpha} \right]\notag\\
 &=P\left[\sup\limits_{x \in (\emph{l},\emph{u})}
     \frac{|(\boldsymbol{B}\mathbf{x})'\boldsymbol{N}|}{\lVert\boldsymbol{B}\mathbf{x}\rVert} < c_{\alpha} \right],
\end{align}
where $\boldsymbol{N} = \boldsymbol{B}^{-1}(\boldsymbol{\beta} - \hat{\boldsymbol{\beta}}) \sim N_2(\boldsymbol{0},\boldsymbol{I})$.

Define $T(\mathbf{x}) = \boldsymbol{B}\mathbf{x}$. Then $T: \mathbb{R}^2 \rightarrow \mathbb{R}^2$ is a linear transformation that maps the line segment in Figure 1 onto the line segment in Figure 2. Therefore, restricting $\mathbf{x}$ to be on the line segment in Figure~1 is equivalent to restricting $\boldsymbol{B}\mathbf{x}$ to lie on the line segment in Figure~2.  Furthermore, we can form a cone, also pictured in Figure 2,  by adding to the line segment all rays through origin that intersect it. Notice that restricting $\boldsymbol{B}\mathbf{x}$ to be on the line is equivalent to restricting it to lie in the cone. This equivalence follows from the fact that the quantity $\frac{|(\boldsymbol{B}\mathbf{x})'\boldsymbol{N}|}{\lVert\boldsymbol{B}\mathbf{x}\rVert}$ in Equation~(\ref{eq:add}) is constant on rays through origin.

Now let $\phi, \, 0 \leq \phi \leq \pi$, be the angle between the vectors $\boldsymbol{B}(1, \,\, \emph{l})'$ and $\boldsymbol{B}(1, \,\, \emph{u})'$, and let $\mathbf{x}_a$ be a vector such that $\boldsymbol{B}\mathbf{x}_a$ cuts $\phi$ into two equal parts. Then the cone can be expressed as: $\{\mathbf{x}: |\rho(\mathbf{x}, \mathbf{x}_a)| \geq \text{cos}(\phi/2)\}$. Clearly this cone has the form of the region $\mathscr{X}$ if we let $\boldsymbol{z}$ = $\mathbf{x}_a$ and $a = \text{cos}(\phi/2)$. Therefore, the interval restriction on the predictor variable in the simple logistic regression model is also a special case of the region $\mathscr{X}$.

\vspace{.3cm}
\emph{\textbf{Remark 3}}. This is given as the corollary in Seppanen $\&$ Uusipaikka \cite{SU92}, which states that the critical value $c_{\alpha}$ is a decreasing function of $a$. This is because for any given numbers $a_1 < a_2$, the region defined by $a_2$ is contained in the region defined by $a_1$, that is, $\mathscr{X}_2 = \{\mathbf{x}: \rho(\mathbf{x}, E)
\geq a_2\} \subseteq \mathscr{X}_1 = \{\mathbf{x}: \rho(\mathbf{x}, E)
\geq a_1\}$.

\subsection{Proof of Theorem 1}
\begin{proof}[Proof of Theorem 1]
Without restriction the $p \times r$ matrix $\boldsymbol{Z}$ can be chosen such that it satisfies $\boldsymbol{Z'F^{-1}Z} = \boldsymbol{I}_r$. Then, there exists a $p \times (p-r)$ matrix $\boldsymbol{D}$  such that $\boldsymbol{D'F^{-1}D} = \boldsymbol{I}_{p-r},$ and $\boldsymbol{Z'F^{-1}D} = \boldsymbol{0}$. That is, the columns of  $\boldsymbol{Z}$ and  $\boldsymbol{D}$ together form an orthonormal basis for the $p$-dimensional space.

Denote $\boldsymbol{z} = \boldsymbol{Z'F^{-1}}\mathbf{x}$ and $\boldsymbol{d} = \boldsymbol{D'F^{-1}}\mathbf{x}$ for every $\mathbf{x} \in \mathbb{R}^p$. Then,
\[
\mathbf{x}'\boldsymbol{F^{-1}ZZ'F^{-1}}\mathbf{x} = \boldsymbol{z}'\boldsymbol{z} = \|\boldsymbol{z}\|^2 \qquad \text{and} \qquad \mathbf{x}'\boldsymbol{F}^{-1}\mathbf{x} = \|\boldsymbol{z}\|^2 + \|\boldsymbol{d}\|^2.
\]
Therefore $\rho(\mathbf{x},E) = [\|\boldsymbol{z}\|^2 / (\|\boldsymbol{z}\|^2 + \|\boldsymbol{d}\|^2)]^{1/2}$, and hence $\mathbf{x} \in \mathscr{X} = \{\mathbf{x}: \rho(\mathbf{x},E) \geq a\}$ if and only if $(\boldsymbol{z}, \boldsymbol{d})' \in \mathscr{D} = \{(\boldsymbol{z}, \boldsymbol{d})': \|\boldsymbol{z}\|^2 / (\|\boldsymbol{z}\|^2 + \|\boldsymbol{d}\|^2)) \geq a^2\}$.

Let $\boldsymbol{u}$ and $\boldsymbol{v}$ be the first $r$ and last $p-r$ components of the vector $(\boldsymbol{Z}, \boldsymbol{D})^{-1}\boldsymbol{F}(\boldsymbol{\hat{\beta}} - \boldsymbol{\beta})$. Then, we have
\begin{align}\label{eq:2.2}
&\{\boldsymbol{\hat{\beta}}: \mathbf{x}'\boldsymbol{\beta} \in \mathbf{x}'\boldsymbol{\hat{\beta}} \pm c_{\alpha}(\mathbf{x}'\boldsymbol{F}^{-1}\mathbf{x})^{1/2},\,\, \forall \mathbf{x} \in \mathscr{X}\}\notag\\
&= \{\boldsymbol{\hat{\beta}}: |\mathbf{x}'(\boldsymbol{\hat{\beta}} - \boldsymbol{\beta})| \leq
                  c_{\alpha}(\mathbf{x}'\boldsymbol{F}^{-1}\mathbf{x})^{1/2},\,\, \forall \mathbf{x} \in \mathscr{X}\}\notag\\
&= \bigg\{\bigg( \begin{array}{c} \boldsymbol{u}\\ \boldsymbol{v} \end{array} \bigg): \bigg|\bigg( \begin{array}{c} \boldsymbol{z}\\ \boldsymbol{d} \end{array} \bigg)'
        \bigg( \begin{array}{c} \boldsymbol{u}\\ \boldsymbol{v} \end{array} \bigg)\bigg| \leq
        c_{\alpha}\bigg\|\bigg( \begin{array}{c} \boldsymbol{z}\\ \boldsymbol{d} \end{array} \bigg)\bigg\|, \quad
        \forall \,\, \bigg( \begin{array}{c} \boldsymbol{z}\\ \boldsymbol{d} \end{array} \bigg) \in \mathscr{D}
 \bigg\}\notag\\
&= \{(\boldsymbol{u}, \boldsymbol{v})': |\boldsymbol{z}'\boldsymbol{u} + \boldsymbol{d}'\boldsymbol{v}| \leq
   c_{\alpha}(\|\boldsymbol{z}\|^2 + \|\boldsymbol{d}\|^2)^{1/2}, \forall \,\, (\boldsymbol{z}, \boldsymbol{d})' \in \mathscr{D}\}.
\end{align}

The joint asymptotic distribution of the vectors $\boldsymbol{u}$ and $\boldsymbol{v}$ is $\boldsymbol{N}_p(\boldsymbol{0},\boldsymbol{I}_p)$. This is because $\boldsymbol{\hat{\beta}} $ is  distributed as $\boldsymbol{N}_p(\boldsymbol{\beta},\boldsymbol{F}^{-1})$ asymptotically, and
\begin{align*}
\mathrm{E}[(\boldsymbol{u}, \boldsymbol{v})'] &= \mathrm{E}[(\boldsymbol{Z}, \boldsymbol{D})^{-1}\boldsymbol{F}(\boldsymbol{\hat{\beta}} - \boldsymbol{\beta})] = (\boldsymbol{Z}, \boldsymbol{D})^{-1}\boldsymbol{F} \mathrm{E}[\boldsymbol{\hat{\beta}} - \boldsymbol{\beta}] = \boldsymbol{0}_{p\times1},\\
\mathrm{Var}[(\boldsymbol{u}, \boldsymbol{v})']
&= (\boldsymbol{Z}, \boldsymbol{D})^{-1}\boldsymbol{F} \boldsymbol{F}^{-1} \boldsymbol{F}'((\boldsymbol{Z}, \boldsymbol{D})^{-1})'\\
&=[(\boldsymbol{Z}, \boldsymbol{D})'\boldsymbol{F}^{-1}(\boldsymbol{Z}, \boldsymbol{D})]^{-1}\\
&=\boldsymbol{I}_p.
\end{align*}

Furthermore, since $\boldsymbol{u}$ and $\boldsymbol{v}$ partition the vector $(\boldsymbol{Z}, \boldsymbol{D})^{-1}\boldsymbol{F}(\boldsymbol{\hat{\beta}} - \boldsymbol{\beta})$, we have
\[
\boldsymbol{u} \sim \boldsymbol{N}_r(\boldsymbol{0},\boldsymbol{I}_r), \quad \boldsymbol{v} \sim \boldsymbol{N}_{p-r}(\boldsymbol{0},\boldsymbol{I}_{p-r}),
\]
and $\boldsymbol{u}$ and $\boldsymbol{v}$ are independent.

Define $r^2 = \|\boldsymbol{z}\|^2 / (\|\boldsymbol{z}\|^2 + \|\boldsymbol{d}\|^2),\,\, 0 \leq r \leq 1$,  and $s = (1 - r^2)^{1/2}$. For all $\boldsymbol{z}$ and $\boldsymbol{d}$, we have the following
\begin{align*}
|\boldsymbol{z}'\boldsymbol{u} + \boldsymbol{d}'\boldsymbol{v}| &\leq |\boldsymbol{z}'\boldsymbol{u}| + |\boldsymbol{d}'\boldsymbol{v}|\\
  &\leq \|\boldsymbol{z}\|\|\boldsymbol{u}\| + \|\boldsymbol{d}\|\|\boldsymbol{v}\| \\
  &\leq (r\|\boldsymbol{u}\| + s\|\boldsymbol{v}\|)(\|\boldsymbol{z}\|^2 + \|\boldsymbol{d}\|^2)^{1/2},
\end{align*}
and since the equality holds when
\[
\boldsymbol{z} = r\frac{\boldsymbol{u}}{\|\boldsymbol{u}\|} \quad \text{and} \quad \boldsymbol{d} = s\frac{\boldsymbol{v}}{\|\boldsymbol{v}\|},
\]
the event in Equation~(\ref{eq:2.2}) is equivalent to the following event:
\begin{equation}\label{eq:2.3}
\{(\boldsymbol{u}, \boldsymbol{v})': r\|\boldsymbol{u}\| + s\|\boldsymbol{v}\| \leq
   c_{\alpha}, \,\, \forall \,\, a \leq r \leq 1\}
\end{equation}

Define $r = \cos(\psi),\,\, 0 \leq \psi \leq \cos^{-1}(a)$. Then  $s = \sin(\psi)$, since $r^2 + s^2 = 1$ by the definitions above. We will also define $w = \|\boldsymbol{u}\|^2 + \|\boldsymbol{v}\|^2$, and $u = \|\boldsymbol{u}\|^2 / (\|\boldsymbol{u}\|^2 + \|\boldsymbol{v}\|^2) = \cos^2(\phi)\, ( 0 \leq \phi \leq \pi/2)$. Then, $w$ follows a $\chi^2$ distribution with $p$ degrees of freedom,  $u$ follows a beta($r/2, (p-r)/2$) distribution, and $w$ and $u$ are independent.

Now the event in Equation~(\ref{eq:2.3}) can be further written as
\begin{equation*}
\{\phi: w\cos^2(\phi - \psi) \leq   c_{\alpha}^2, \,\, \forall \,\, \psi \in \Psi = [0,\, \cos^{-1}(a)]\},
\end{equation*}
or equivalently,
\begin{equation}\label{eq:2.4}
\{\phi: w\, \sup_{\psi \in \Psi}\,\cos^2(\phi - \psi) \leq   c_{\alpha}^2\}.
\end{equation}

Define $G = w\, \sup_{\psi \in \Psi}\,\cos^2(\phi - \psi)$. Then $c_{\alpha}^2$ is the upper $\alpha$ point of the distribution function of the random variable $G$. Notice that
\[
\sup_{\psi \in \Psi}\,\cos^2(\phi - \psi)=
\begin{cases}
1, &\text{if} \,\, \phi \in \Psi\\
\cos^2(\phi - \cos^{-1}(a)), &\text{if}\,\, \phi \notin \Psi.
\end{cases}
\]
Then the distribution function of the random variable $G$ is given by
\begin{align}\label{eq:2.5}
P(G \leq g) &= P(w\, \sup_{\psi \in \Psi}\,\cos^2(\phi - \psi) \leq g) \notag\\
            &= P(w \leq g) + P(g < w < g/\sup_{\psi \in \Psi}\,\cos^2(\phi - \psi))\notag\\
            &= P(w \leq g) + P(g < w < g/\cos^2(\phi - \cos^{-1}(a)), \,\, \cos^{-1}(a) < \phi \leq \pi/2)
\end{align}
Notice the event $\{g < w < g/\cos^2(\phi - \cos^{-1}(a))\}$ can be written as $\{\sqrt{g/w} < 1 < \sqrt{g/w} / \cos(\phi - \cos^{-1}(a))\}$, which implies $\phi > \cos^{-1}(a) + \delta$,  where $\delta$ is defined as $\delta = \cos^{-1}(\sqrt{g/w})$. Hence Equation~(\ref{eq:2.5}) can be further written as
\begin{equation}\label{eq:2.6}
P(w \leq g) + P(g < w < g/\cos^2(\pi/2 - \cos^{-1}(a)), \,\, \cos^{-1}(a) + \delta < \phi \leq \pi/2).
\end{equation}
Furthermore, the event $\{\cos^{-1}(a) \, + \, \delta < \phi \leq \pi/2\}$ is equivalent to the event $\{0 \leq \cos^2(\phi) < \cos^2(\cos^{-1}(a) + \delta)\}$, and
\[
\cos^2(\cos^{-1}(a) + \delta)=\{a\sqrt{g/w} - [(1 - a^2)(1 - (\sqrt{g/w})^2)]^{1/2} \}^2.
\]
Define a function $m(t) = \{at -  [(1 - a^2)(1 - t^2)]^{1/2}\}^2$. Then
\[
\{\cos^{-1}(a) \, + \, \delta < \phi \leq \pi/2\} = \{0 \leq u \leq m(\sqrt{g/w})\}
\]
Also $\cos^2(\pi/2 - \cos^{-1}(a))$ in Equation~(\ref{eq:2.6}) can be simplified as $1 - a^2$. Therefore Equation~(\ref{eq:2.6}) can be written as the following, which gives the distribution function of $G$,
\begin{align*}
P(G \leq g)&=
P(w \leq g) + P(g < w < g/(1 - a^2), \,\, 0 \leq u \leq m(\sqrt{g/w}))\\
&= F(g) + \int_g^{g/(1-a^2)}H(m(\sqrt{g/w}))f(w)dw,
\end{align*}
where $H$ is the distribution function of the beta($r/2, (p-r)/2$) distribution. $F$ and $f$ are the distribution and density functions of the $\chi^2_p$ distribution, respectively.
\end{proof}

\subsection{Conservative Two-sided Bands Over Convex Regions}
Constraint sets for predictor variables are not necessarily of the form~(\ref{eq:1.5}). For example, it is often the case that the practical range of values for each predictor variable are specified by using a lower and  upper bound, so the resulting constraint region is a rectangular region. The following theorem, which is due to  Seppanen $\&$ Uusipaikka \cite{SU92}, provides a conservative confidence bands for such convex sets. The theorem is restated here and the proof can be found in Seppanen $\&$ Uusipaikka \cite{SU92}. Let $\mathscr{X}^*$ be a given subset of the $p$-dimensional space, not of the form~(\ref{eq:1.5}).

\begin{theorem}
 If $\mathscr{X}^*$ is a convex set generated by the vectors $\mathbf{x}_1, \mathbf{x}_2, \ldots, \mathbf{x}_k$ ($p \leq k$) and  $\mathbf{x}_0$ is a given vector of $\mathscr{X}^*$ such that $a = \min_{1 \leq i \leq k}\rho(\mathbf{x}_i, \mathbf{x}_0) > 0$, then the confidence bands~\eqref{eq:1.8} and \eqref{eq:1.10} are conservative $(1 - \alpha)$-level confidence bands over the region $\mathscr{X}^*$, if $c_{\alpha}$ is equal to the square root of the upper $\alpha$ point of the distribution function~\eqref{eq:2.1}.

\end{theorem}
According to Remark 3 in section~2.1, we should choose the vector $\mathbf{x}_0$ in Theorem~2 such that $a = \inf\{\rho(\mathbf{x}, \mathbf{x}_0): \mathbf{x} \in \mathscr{X}^*\}$ is as large as possible. Finding the ``best'' vector $\mathbf{x}_0$ which yields the largest possible value of $a$ requires an iterative process. However, some existing computer programs, such as $\boldsymbol{R}$, which will be used in the later examples,  greatly simplify the process and allow us to find the ``best'' vector $\mathbf{x}_0$  very fast. A sample of $\boldsymbol{R}$ code for finding the ``best'' choice of $\mathbf{x}_0$ is included in appendix.

\section{Examples}
\subsection{Simple Linear Logistic Regression}
We will use the genetic toxicity data obtained by LaVelle \cite{L86}. These were the same data that were considered by Piegorsch $\&$ Casella \cite{PC88} and Kerns \cite{K15}. They are used here again to illustrate the proposed method, and to compare the new method with the existing methods provided by  Piegorsch $\&$ Casella and Kerns. The study in LaValle \cite{L86} investigated frameshift mutagenesis in bacterial assays. Here, the bacterium under study is the \emph{E}. \emph{coli}, strain 343/435, and the findings for a control and five doses of the suspected mutagen, 9-Aminoacridine (9-AA) are reported in Table~\ref{ta:3.1}.
\begin{table}
\tbl{Mutagenicity of 9-Aminoacridine in \emph{E}. \emph{coli} strains 343/435}
{\begin{tabular}[l]{c|c|c|c|c|c|c}
\hline
Dose &-$^{\rm a}$       &.8   &2.4   &8.0      &24    &80\\\hline
Log-dose      &-1.374    &-.223  &0.875 &2.079   &3.178  &4.382\\\hline
Response       &7/96     &28/96  &64/96 &54/96  &81/96 &96/96\\\hline
\end{tabular}}
\tabnote{$^{\rm a}$The first data pair corresponds to a zero-dose control. The log-dose for this datum was calculated using consecutive-dose average spacing (Margolin \emph{et al.}, 1986)}
\label{ta:3.1}
\end{table}

As mentioned in Piegorsch $\&$ Casella's paper, a simple linear logistic regression provides a good model for the data.  So, we fit a simple linear logistic model to the data using the log-dose level as our predictor variable $x$, and the ML estimates from the logistic fit
are $\hat{\beta}_0 = -.789$, and $\hat{\beta}_1 = .854$. The inverse of the Fisher
information matrix is
$\boldsymbol{F^{-1}}=\left[ \begin{array}{cc} 0.017 &-0.005 \\ -0.005
    & 0.005 \end{array}\right]$, and the matrix $\boldsymbol{B}$, defined in Remark 2, is $\left[ \begin{array}{cc} 0.128 &-0.027 \\ -0.027
    & 0.063 \end{array}\right]$.

The computer language used to fit the simple logistic regression is $\boldsymbol{R}$, version 3.2.1. If we construct confidence bands over the whole number line, the bands will be unnecessarily wide. So it is often of interest to restrict the predictor variable to a certain interval, and then construct narrower confidence bands over the interval. Also it is often
noted that human exposure to environmental toxins usually occurs at
low dose levels.  Therefore, we will focus on the constrained intervals that are concentrated toward the lower end of possible values, and build confidence bands over the selected intervals.

Three constrained intervals were studied by Piegorsch $\&$ Casella \cite{PC88} and Kerns \cite{K15}, and they will be restudied here. For example, consider one of the intervals, $(-1.3,\,.8)$.  It is straightforward to calculate the angle between the vectors $\boldsymbol{B}(1, \,\, -1.3)'$ and $\boldsymbol{B}(1, \,\, .8)'$, and then the value of $a$. For this interval, the angle and $a$ are found to be $\phi = .809$ and $a = .9193$ , which yields a critical value of  $c_{\alpha} = 2.206$ for the bands at a $95\%$ confidence level. Therefore, a $95\%$ two-sided confidence band for $p(\mathbf{x})$ when $x \in (-1.3, \, .8)$ is given by

\begin{equation*}
 p(\mathbf{x}) \in
 \{1 + \text{exp}[(.789-.854x) \pm 2.206(\mathbf{x}'\boldsymbol{F}^{-1}\mathbf{x})^{1/2}]\} ^{-1}.
\end{equation*}

 The values of $a$ and critical points $c_{\alpha}$ based on the proposed method are reported in Table~\ref{ta:3.2} for three selected intervals, along with the
critical points given by Piegorsch $\&$ Casella \cite{PC88} and Kerns \cite{K15}. Values of the critical
point, based on Scheff\'{e}'s method, when there are no restrictions on
$x$, are also given here. The constraint considered in this example has the form $(\emph{l},\,\,u\emph{})$ with $\emph{l}$ and $\emph{u}$ being given constants, which is the restricted interval considered by Kerns~\cite{K15} in the simple logistic regression case.  It can be seen from the table that the results based on the proposed method are consistent with the results using Kerns' method. This should be expected because as stated in Remark 2, the interval restriction on $x$ in the simple logistic regression model is a special case of the region $\mathscr{X}$. It is also clear from the table that the proposed method is able to produce smaller critical values and hence improves the method given by Piegorsch $\&$ Casella in this case.

\begin{table}
\tbl{Critical Values for Two-sided Bands with $\alpha=.05$}
{\begin{tabular}{c c c c c c c}
\hline
\multirow{2}{*}{Restricted Intervals}  &Kerns  & &\multicolumn{2}{c}{Proposed method} & & Piegorsch and Casella  \\\cline{2-2}\cline{4-5}\cline{7-7}
                                       &$c_{\alpha}$  &  &$a$   &$c_{\alpha}$       &             &$c_{\alpha}$ \\\hline
$(-\infty, \,\, \infty)$               &2.447     &     & 0        &2.447  &   &2.447\\\hline
(-1.3,\,\,2.0)                         &2.344     &     &.7233    &2.344  &   &2.445 \\\hline
(-1.3,\,\,0.8)                         &2.206     &     &.9193    &2.206  &   &2.274 \\\hline
(-1.3,\,\,-.2)                         &2.067     &    &.9887    &2.067   &   &2.170 \\\hline
\end{tabular}}
\label{ta:3.2}
\end{table}

\subsection{Multiple Linear Logistic Regression}
For the multiple logistic regression case, we will use the data set, called ICU,  given by Hosmer, Lemeshow and Sturdivant \cite{HL13}. The data set consists of a sample of 200 subjects who were part of a much larger study on survival rates following admission to an adult intensive care unit (ICU). The outcome is vital status, alive or dead, coded as $0/1$ respectively, under the variable name ``STA''. 20 independent variables were investigated in the study, which include age, sex, race, etc.  The question of interest is to determine how the probability of survival after being discharged is related to these independent variables, and a major goal of the study is to develop a multiple logistic regression model to investigate the effect of these independent variables on the survival probability of  the patients after being discharged.  For this illustrative example, we will investigate the effect of the two variables, age and systolic blood pressure at ICU admission, which were under the variable names ``AGE" and ``SYS", respectively.

The sample size in the example is $n=200$, therefore it is sufficiently
large for the asymptotic approximation to hold, as the Monte
Carlo simulation studies in Section~4 will show.  The software $\boldsymbol{R}$ is  used again to fit the multiple logistic regression model using $y$ (vital status) as the response variable, $x_1$ (age) and $x_2$ (blood pressure) as the independent variables. The ML estimates from the multiple logistic fit
are $\hat{\beta}_0 = -.962$, \,\, $\hat{\beta}_1 = .028$, and $\hat{\beta}_2 = -.017$. Both age and blood pressure are found to have a significant effect on the vital status, with \emph{p}-values of $.00838$ and $.00407$ respectively. The inverse of the Fisher
information matrix is
$\boldsymbol{F^{-1}}=\left[ \begin{array}{ccc} 1.001 &-.0072 &-.0041 \\
                                               -.0072 &.00012 &.000001\\
                                               -.0041 &.000001 &.00003
     \end{array}\right]$.

In the data, the variable $x_1$, age, ranges between 16 and 92 years, and the variable $x_2$, blood pressure, ranges between 36 and 256. If we restrict $x_1$ to the interval ($16,\,\,92$), and $x_2$ to the interval ($36,\,\,256$), then the constraint is a rectangular region, which is a convex set, but not of the form $\mathscr{X}$. We apply the iterative procedure in Theorem~2 to this example to search for the optimal $\mathbf{x}_0$. The $\boldsymbol{R}$ code for finding the best choice of $\mathbf{x}_0$ and the corresponding value for $a$ is given in appendix.  Our $\boldsymbol{R}$ program calculates $a = .2383$,  which occurs at $\mathbf{x}_0 = (1, \,\, 62.91, \,\, 124.62)'$. The resulting critical value is found to be $c_{\alpha} = 2.789$ at a $95\%$ confidence level. Therefore, a $95\%$ conservative confidence band for
$p(\mathbf{x})$ when $\mathbf{x}$ is constrained in the specified region is given by
\begin{equation*}
 p(\mathbf{x}) \in
 \{1 + \text{exp}[(-.962 + .028x_1 - .017x_2) \pm 2.789(\mathbf{x}'\boldsymbol{F}^{-1}\mathbf{x})^{1/2}]\} ^{-1}.
\end{equation*}

A patient is considered as having high blood pressure if his/her systolic blood pressure falls between 140 and 160, and is considered as having hypertensive crisis if his/her blood pressure is higher than 180. Patients with systolic blood pressure falling below 120 are considered as normal. Furthermore, we would like to investigate the difference between the young patients and old patients. Several regions are studied based on these considerations, and corresponding values of $a$ and $c_{\alpha}$ are given in Table~\ref{ta:3.3}, along with the critical value using Scheff\'{e}'s method when there are no restrictions on the independent variables.  It can be seen from the table
that there is clear improvement in the width of confidence bands when the independent variables are restricted to smaller regions. For example, one particular region considered is $x_1 \in (20,\,\, 40)$ and $x_2 \in (140,\,\,160)$, in which case we focus our attention on patients who are young and have high blood pressure. The resulting critical value is 2.220, and hence a saving in bands width as great as 21$\%$ over the unrestricted critical value given by Scheff\'{e}'s method has been achieved.

\begin{table}
\tbl{Critical Values for Two-sided Bands \\with $\alpha=.05$}
{\begin{tabular}{c c c c}
\hline
Restriction on $x_1$      &Restriction on $x_2$         &$a$       &$c_{\alpha}$  \\\hline
$(-\infty, \,\, \infty)$  &$(-\infty, \,\, \infty)$     &0         &2.795\\\hline
(16,\,\,92)               &(36,\,\,256)                 &.2383     &2.789 \\\hline
(20,\,\,40)               &(140,\,\,160)                &.9731     &2.220 \\\hline
(50,\,\,80)               &(140,\,\,160)                &.7917    &2.557 \\\hline
(20,\,\,40)               &(30,\,\,120)                &.8658     &2.468 \\\hline
(50,\,\,80)               &(30,\,\,120)                &.7007     &2.634 \\\hline
(20,\,\,40)               &(180,\,\,250)                &.9560     &2.283 \\\hline
(50,\,\,80)               &(180,\,\,250)                &.9200     &2.374 \\\hline
\end{tabular}}
\label{ta:3.3}
\end{table}

\section{Monte Carlo Simulation}
The confidence bands in Equations~(\ref{eq:1.8}) and
(\ref{eq:1.10}) are constructed based on the asymptotic properties of
the maximum likelihood estimators in the multiple logistic model
assuming the sample size is large, therefore an assessment of the small sample performance is needed. We conducted Monte Carlo
simulation studies in $\boldsymbol{R}$ to estimate the actual coverage level of the bands when the sample size is small. The simulation study was performed for the case of one predictor variable. Five different values of $\boldsymbol{\beta}$ were chosen to evaluate the performance of the proposed method:
$[-2,\, .3]',\, [0,\, 1.5]',\,[2,\, 5]',\,[-.2,\,
-.3]',\, [-2, -4]'$.
These selected values for $\boldsymbol{\beta}$ were the same values considered by Kerns \cite{K15}, and they represent five different forms of the probability response
functions: slowly increasing, moderately increasing, fast increasing,
slowly decreasing, and fast decreasing, respectively.

Three different intervals on the independent variable were examined: narrow, wide, and extremely wide ("unrestricted"). The endpoints for each interval were obtained  by inverting Equation~(\ref{eq:1.1}), that is, they were calculated based on the formula $x=(\text{log}_e[p/(1-p)]-\beta_0)/\beta_1$.  For a given $\boldsymbol{\beta}$, values of $p$ were selected to produce the endpoints of these intervals. In particular, $p=.3,\, .7$  were used for computing the endpoints of the narrow intervals, $p=.1,\, .9$ were selected for the wide intervals, and $p=10^{-10},\, 1-10^{-10}$ were chosen for the very wide intervals. These were also the same intervals studied by Kerns \cite{K15}, and for completeness, they are presented here in Table~\ref{ta:4.1}.

\begin{table}
 \tbl{Restricted Intervals For Monte Carlo Simulation}
 {\begin{tabular}{cccc}
 \hline\noalign{\smallskip}
   $\boldsymbol{\beta}$   &Narrow Interval &Wide Interval  &\textquotedblleft Unrestricted Interval\textquotedblright \\\hline
   $[-2,\, .3]'$   &$(3.842, \, 9.491)$   &$(-.657, \, 13.991)$    &$(-70.086, \, 83.420)$\\
   $[0,\, 1.5]'$   &$(-.565, \, .565)$   &$(-1.465, \, 1.465)$    &$(-15.351, \, 15.351)$\\
   $[2,\, 5]'$   &$(-.569, \, -.231)$   &$(-.839, \, -.039)$    &$(-5.005, \, 4.205)$\\
   $[-.2,\, -.3]'$   &$(-3.491, \, 2.158)$   &$(-7.991, \, 6.657)$    &$(-77.420, \, 76.086)$\\
   $[-2,\, -4]'$   &$(-.712, \, -.288)$   &$(-1.049, \, .049)$    &$(-6.256, \, 5.256)$\\\hline
    \end{tabular}}
 \label{ta:4.1}
 \end{table}

Once the intervals were set, equidistant values bounded in the specified ranges were generated as the values for the
predictor variable $x$ at four different sample sizes ($n = 25, 50, 100, 150$). The
dichotomous response variable $Y$ was generated based on the
following: First a uniform
$(0,1)$ random variable with the specified sample size was simulated, and then $Y$ was determined from the probability of success, $p(x) = 1/(1+\text{exp}(\beta_0+\beta_1x))$: $Y=1$ if the uniform random variable was less than
$p(x)$; $Y=0$ otherwise.

For each simulated data set, the ML estimate of the parameter $\boldsymbol{\beta}$ and the inverse of the Fisher information matrix $\boldsymbol{F}^{-1}$ were obtained and used to form confidence bands based on the proposed method to evaluate the coverage probability.

The number of iterations in each Monte Carlo simulation was chosen as $N=5000$, and  three nominal error rates $\alpha$ ($.01, .05, .10$) were considered in the study. The Monte Carlo error,  $1 - $ the
 estimated coverage probability, was estimated and the results are presented in Table~\ref{ta:4.2}. It is clear from the table that the proposed procedure is conservative for
 small samples, but the error approaches the nominal level as the
 sample size increases. Generally, the error reaches the nominal level
 when the sample size is 100, but in some cases it could be as small
 as 50. It can be seen from the table that there is no noticeable difference in errors
 between three nominal error rates. It is also noted that three different intervals exhibit similar coverage results.

\begin{table}
 \tbl{Estimated Monte Carlo Errors For the Confidence Bands}
 {\begin{tabular}{c|c|ccc|ccc|ccc}
 \hline\noalign{\smallskip}
   \multirow{2}{*}{$\boldsymbol{\beta}$}
     &\multirow{2}{*}{Sample size $n$}
     &\multicolumn {3}{c|}{Narrow}  &\multicolumn {3}{c|}{Wide}
     &\multicolumn {3}{c}{\textquotedblleft Unrestricted \textquotedblright} \\\cline{3-11}

    & &$\alpha=.01$ &$.05$ &$.10$ &$.01$ &$.05$ &$.10$ &$.01$ &$.05$ &$.10$ \\\hline
  \multirow{4}{*}{$[-2,\, .3]'$}    &25 &$.001$ &$.022$ &$.060$  &$.001$ &$.024$ &$.061$  &$.002$ &$.025$ &$.052$\\
                                    &50 &$.005$ &$.038$ &$.076$  &$.004$ &$.036$ &$.081$   &$.007$ &$.037$ &$.066$\\
                                    &100 &$.008$ &$.045$ &$.086$ &$.008$  &$.039$ &$.088$  &$.009$  &$.039$ &$.080$\\
                                    &150 &$.009$ &$.049$ &$.096$ &$.008$  &$.051$ &$.096$  &$.010$  &$.043$ &$.081$\\\hline
  \multirow{4}{*}{$[0,\, 1.5]'$}    &25 &$.001$ &$.022$ &$.065$  &$.004$ &$.034$ &$.064$  &$.002$ &$.030$ &$.065$\\
                                    &50 &$.005$ &$.040$ &$.076$  &$.007$ &$.035$ &$.078$  &$.005$ &$.036$ &$.066$\\
        &100 &$.006$ &$.043$ &$.086$ &$.009$ &$.040$ &$.092$  &$.008$ &$.044$ &$.080$\\
        &150 &$.008$ &$.042$ &$.097$ &$.008$ &$.042$ &$.092$  &$.011$ &$.047$ &$.081$\\\hline
  \multirow{4}{*}{$[2,\, 5]'$}    &25 &$.001$ &$.022$ &$.063$ &$.006$ &$.033$ &$.063$  &$.004$ &$.022$ &$.062$\\
        &50 &$.005$ &$.035$ &$.089$  &$.007$ &$.033$ &$.078$   &$.006$ &$.036$ &$.065$\\
        &100 &$.008$ &$.045$ &$.090$ &$.007$ &$.046$ &$.094$   &$.010$ &$.039$ &$.078$\\
        &150 &$.007$ &$.046$ &$.099$  &$.008$ &$.045$ &$.093$  &$.011$ &$.044$ &$.092$\\\hline
  \multirow{4}{*}{$[-.3,\, -.2]'$}    &25 &$.001$ &$.023$ &$.069$  &$.005$ &$.029$ &$.061$  &$.005$ &$.018$ &$.058$\\
        &50 &$.005$ &$.036$ &$.080$  &$.007$ &$.031$ &$.074$  &$.007$ &$.039$ &$.068$\\
        &100 &$.008$ &$.044$ &$.092$ &$.008$ &$.048$ &$.096$  &$.009$ &$.042$ &$.090$\\
        &150 &$.008$ &$.045$ &$.099$ &$.009$ &$.046$ &$.089$  &$.010$ &$.046$ &$.094$\\\hline
  \multirow{4}{*}{$[-4,\, -2]'$}    &25 &$.002$ &$.021$ &$.062$  &$.007$ &$.033$ &$.064$  &$.007$ &$.024$ &$.063$\\
        &50 &$.004$ &$.032$ &$.073$ &$.008$ &$.035$ &$.074$   &$.008$ &$.034$ &$.074$\\
        &100 &$.009$ &$.044$ &$.090$  &$.009$ &$.041$ &$.091$  &$.008$  &$.043$ &$.088$\\
        &150 &$.010$ &$.045$ &$.094$  &$.012$ &$.048$ &$.088$  &$.011$ &$.047$ &$.092$\\\hline
    \end{tabular}}
 \label{ta:4.2}
 \end{table}

\appendices
\section{$R$ code for finding the best vector $\mathbf{x}_0$ and the value of $a$ in Example 3.2.}

\begin{lstlisting}
#Set up the data set
icu_data=read.csv("C:/Users/xlu/Dropbox/on-going research
         /Logistic_bands_multiple/icu_example2.csv",header=T)
mydata=data.frame(sta=icu_data$STA, age=icu_data$AGE,
                  blood_pressure=icu_data$SYS, rate=icu_data$HRA)

#Run logistic regression on the data
output <- glm(sta~age+blood_pressure,data=mydata,family=binomial)

#Extract the inverse of the Fisher Information matrix
Finv <- summary(output)$cov.scaled

#Set the interval on the explanatory variable x1
lower1=16
upper1=92

#Set the interval on the explanatory variable x2
lower2=36
upper2=256

#Define the vectors that generate the restriction region
x1=c(1,lower1,lower2)
x2=c(1,lower1,upper2)
x3=c(1,upper1,lower2)
x4=c(1,upper1,upper2)
x=cbind(x1,x2,x3,x4)


#Initiating
N1=500
N2=500
minrho=matrix(0,N1,N2)
h1=seq(lower1,upper1, length.out=N1)
h2=seq(lower2,upper2, length.out=N2)

#Run the iterative process to find the value of x_0 which
#yields the maximum value of rho.

for (i in 1:N1)
{
  for (k in 1:N2)
  {
    x0=c(1,h1[i],h2[k])
    E00=t(x0) %*% Finv %*% x0
    E=array(0,4)
    rho=array(0,4)

    #Calculate the value of rho for a given x_0
    for (j in 1:4)
    {
      E[j]=t(x[,j])%*% Finv %*% x[,j]
      E10=t(x[,j])%*% Finv %*% x0
      rho[j]=E10/sqrt(E[j]*E00)
    }
    #Store the value of rho in the matrix "minrho"
    minrho[i,k] <- ifelse(min(rho)>0,min(rho),0)
  }
}

#The maximum value of rho is the value of a
a=max(minrho)

#Locate the vector of x_0 that yields the maximum of rho
location <- which(minrho==max(minrho), arr.ind=TRUE)
x0_best=c(1,h1[location[1]],h2[location[2]])

#Display the results
cat("a =", a, "which occurs at x_0 =", "[",x0_best,"]","\n")

\end{lstlisting}

\end{document}